\documentclass{ws-procs}
\usepackage{amsfonts,amsmath,amssymb}

\begin{document}


\title{Exotic bialgebras : non-deformation quantum~groups}

\author{D. Arnaudon}
\address{Laboratoire d'Annecy-le-Vieux de Physique Th{\'e}orique LAPTH,\\
UMR 5108 du CNRS associ\'ee \`a l'Universit\'e de Savoie,\\ 
  BP 110, F-74941 Annecy-le-Vieux~Cedex, France,\\
  daniel.arnaudon@lapp.in2p3.fr}  

\author{A. Chakrabarti}
\address{Centre de Physique Th\'eorique,
  \'Ecole Polytechnique,\\ F-91128 Palaiseau Cedex, France,\\ 
  chakra@cpht.polytechnique.fr}  

\author{V.K. Dobrev}
\address{School of Informatics, 
University of Northumbria,\\ Newcastle-upon-Tyne NE1 8ST, UK,\\
vladimir.dobrev@unn.ac.uk,\\ 
Permanent address:\\ 
Institute of Nuclear Research and Nuclear Energy,\\ 
Bulgarian Academy of Sciences,\\
  72 Tsarigradsko Chaussee, 1784 Sofia, Bulgaria,\\ 
  dobrev@inrne.bas.bg}

\author{S.G. Mihov}
\address{Institute of
  Nuclear Research and Nuclear Energy,\\  
  Bulgarian Academy of Sciences,\\ 
  72 Tsarigradsko Chaussee, 1784 Sofia, Bulgaria,\\ 
  smikhov@inrne.bas.bg}  


\maketitle

\abstracts{
In the classification of solutions of the Yang--Baxter equation, there are
solutions that are not deformations of the trivial solution (essentially
the identity). We consider the algebras defined by these 
solutions, and the corresponding dual algebras. We then study the 
representations of the latter. We are also interested in the Baxterisation
of these $R$-matrices and in the corresponding quantum planes.
}


\def\cal#1{{\mathcal #1}}

\newcommand{\II}{{\mathbb I}}
\def\CC{{\mathbb C}}
\def\NN{{\mathbb N}}
\def\QQ{{\mathbb Q}}
\def\RR{{\mathbb R}}
\def\ZZ{{\mathbb Z}}
\def\cA{{\cal A}}          \def\cB{{\cal B}}          \def\cC{{\cal C}}
\def\cD{{\cal D}}          \def\cE{{\cal E}}          \def\cF{{\cal F}}
\def\cG{{\cal G}}          \def\cH{{\cal H}}          \def\cI{{\cal I}}
\def\cJ{{\cal J}}          \def\cK{{\cal K}}          \def\cL{{\cal L}} 
\def\cM{{\cal M}}          \def\cN{{\cal N}}          \def\cO{{\cal O}}
\def\cP{{\cal P}}          \def\cQ{{\cal Q}}          \def\cR{{\cal R}} 
\def\cS{{\cal S}}          \def\cT{{\cal T}}          \def\cU{{\cal U}}
\def\cV{{\cal V}}          \def\cW{{\cal W}}          \def\cX{{\cal X}}
\def\cY{{\cal Y}}          \def\cZ{{\cal Z}}
\def\bA{{\bar A}}          \def\bB{{\bar B}}          \def\bC{{\bar C}}
\def\bD{{\bar D}}          \def\bE{{\bar E}}          \def\bF{{\bar F}}
\def\bG{{\bar G}}          \def\bH{{\bar H}}          \def\bI{{\bar I}}
\def\bJ{{\bar J}}          \def\bK{{\bar K}}          \def\bL{{\bar L}} 
\def\bM{{\bar M}}          \def\bN{{\bar N}}          \def\bO{{\bar O}}
\def\bP{{\bar P}}          \def\bQ{{\bar Q}}          \def\bR{{\bar R}} 
\def\bS{{\bar S}}          \def\bT{{\bar T}}          \def\bU{{\bar U}}
\def\bV{{\bar V}}          \def\bW{{\bar W}}          \def\bX{{\bar X}}
\def\bY{{\bar Y}}          \def\bZ{{\bar Z}}
\def\tL{{\tilde L}}
\def\diag{\mathop{\rm diag}\nolimits}
\def\qmbox#1{\qquad\mbox{#1}\quad}

\newcommand{\un}{\mbox{1\hspace{-1mm}I}}
\newcommand{\eps}{{\varepsilon}}
\newcommand{\id}{\mbox{id}}
\newcommand{\Rtriang}[9]{
  {\left(
    \begin{array}{cccc}
      1 & \ #1\  & \ #2\  & \ #3 \ \cr
      & #4 & #5 & \ #6\  \cr
      & & \ #7\  & \ #8\  \cr
      & & & #9
    \end{array}
  \right)
  }}

\def\lg{\left\langle}
\def\rg{\right\rangle}

\def\ta{\tilde a}  \def\td{\tilde d}
\def\tA{\tilde A}  \def\tD{\tilde D}

\def\eqn#1{\begin{equation}\label{#1}}

\newcommand{\eqna}[1]{\begin{subequations} \label{#1}
\begin{eqnarray}}
\def\eena{\end{eqnarray}
\end{subequations}}

\def\hR{\hat R} 


\section{Introduction: {\bf Five}~ 
rank-4 $R$-matrices which are ~{\bf not}~ \\deformations of the identity
}
We are interested in the study of algebraic structures coming from
$R$-matrices (solutions of the Yang--Baxter equation) 
that are not deformations of classical ones (i.e., the
identity up to signs). Those matrices were obtained by Hlavat\'y
\cite{Hlavaty} and
are also in the classification of Hietarinta \cite{Hietarinta}.
There are five such $R$-matrices that are invertible. 
Most of this text is based on the articles \cite{ACDM1, ACDM2, ACDM3}. 

The first cases of interest are
\begin{equation}
  R_{H2,3}
 = \Rtriang{x_1}{x_2}{x_3}
  {1}{0}{x_2}
  {1}{x_1}
  {1}
\end{equation}
denoted by 
\begin{itemize}
\item 
\emph{Exotic 1} (E1) if $x_1=-x_2=-h, h_3 \neq -h^2$ 
\cite{ACDM1}, 
\item 
\emph{Exotic 2} (E2) if $x_1\neq -x_2$ \cite{ACDM1}.
\end{itemize}

Note that if  $x_1=-x_2=-h, h_3 = -h^2$, then this is 
$R_{H1,3}$ with $g=-h$ in the notation of Hietarinta, i.e., Jordanian
deformation with two parameters.) 

The case denoted by \emph{Exotic 3} (E3) is related to the $R$-matrix
\begin{equation}
  R_{S0,2} = \Rtriang{0}{0}{1}
  {-1}{0}{0}
  {-1}{0}
  {1}
\end{equation}
and the algebraic structure was studied in \cite{ACDM1}.

The cases denoted by ``S03'' and ``S14'' come from the $R$-matrices
\begin{equation}
 \label{eq:S03}
 R_{S0,3}\equiv
 \left(
   \begin{array}{rrrr}
     1 & 0 & 0 & 1 \cr
     0 & 1 & 1 & 0 \cr
     0 & 1 & -1 & 0 \cr
     -1 & ~0 & ~0 & ~1 \cr 
   \end{array}
 \right)
 \ ,\qquad
 R_{S1,4}\equiv
 \left(
   \begin{array}{cccc}
     ~0~ & ~0~ & ~0~ & ~q~ \cr
     0 & 0 & 1 & 0 \cr
     0 & 1 & 0 & 0 \cr
     q & 0 & 0 & 0 \cr 
   \end{array}
 \right) \;.
\end{equation}
The structures corresponding to S03 and S14 were studied in \cite{ACDM2}.

\section{
{Algebra and co-algebra structures}
}

The algebra relations are obtained using
\begin{equation}
  R_{12} T_1 T_2 = T_2 T_1 R_{12}
  \label{eq:RTT}
\end{equation}
\begin{equation}
  \mbox{with}
  \qquad
  T = \begin{pmatrix}
    ~a~ & ~b~ \cr c & d
    \end{pmatrix}
\qquad 
T_1 = T \otimes 1
\qquad 
T_2 = 1 \otimes T
\end{equation}
\\
The coalgebra structure is as usual given by
$ \Delta(T) = T \stackrel{.}{\otimes} T $
\\
For $R= \un_2 \otimes \un_2$, the   relations 
$RTT=TTR$ are just the commutativity of $a,b,c,d$ 
($ab=ba$, $ac=ca$, ...),   
i.e., $T$ is a matrix of commuting objects.
For $R$ of ``quantum'' type (with two parameters $q,p$), the relations are
\begin{eqnarray}
  && ab = qba \qquad ac = pca \qquad bd = pdb \\
  && cd = qdc \qquad qbc = pcb \qquad ad - da = (q-p^{-1}) bc 
\end{eqnarray}
which are deformations of simple commutativity relations. The
relations we will obtain in the exotic cases are not such deformations.
\\
For the following
we use the generators:
$ \tilde{a}  = \frac12 (a+d), \ \ \ \ \tilde{d}  = \frac12 
(a-d)
$.

\section{The E1 Case}

\subsection{Algebra and coalgebra relations}

The algebra relations (\ref{eq:RTT}) in the E1 case are:
\begin{alignat}{3}
  &\tilde{a}c  = c\tilde{a} =
  \tilde{d}c  = c\tilde{d} = \tilde{a}\tilde{d} =
  \tilde{d}\tilde{a}  = 0\ ,\quad
 && cb  = bc\ ,\quad
 &&
 c^2  = 0\ , \\
 & \tilde{d}b  = b\tilde{d} + 2h\tilde{d}^2 + hbc \quad
 && \tilde{a}b  = b\tilde{a}, 
 &&
\end{alignat} 
and a basis of the enveloping algebra is
\begin{equation}
 \label{eq:PBWbb}
b^n \ta^k\,, \quad b^n \td^\ell \,, \quad b^n c\,, \qquad 
n,k\in \ZZ_+ \,, \quad \ell\in\NN\ .
\end{equation}
$h_3$ is not explicitly in the relations, but some relations exist
uniquely because $h_3\neq -h_1^2$.   
This algebra has ideals
$I  = {\cA_1} b \tilde d \oplus {\cA_1} \td^2 \oplus {\cA_1} b c $,~~ 
$I_2  = {\cA_1} \td^2 \oplus {\cA_1} b c $ and 
$I_1  = {\cA_1} b c $
such that 
$ I_1 \subset I_2 \subset I \subset \cA_1 $.
The dual algebra, or, more precisely, the \emph{algebra in
duality}, 
is defined via non-degenerate pairing ~$\lg\cdot,\cdot\rg$, 
consistent with the co-algebraic structure, i.e.
such that, 
\eqna{dub} 
&\lg \ u \ , \ ab \ \rg \ =\ \lg \ \Delta(u) \ ,
\ a \otimes b \ \rg \ , \ \ \ \lg \ uv \ , \ a \ \rg \ =\ 
\lg \ u \otimes v \ , \ \Delta (a) \ \rg \qquad \\ 
&\lg 1_\cU , a \rg \ =\ \varepsilon_\cA (a) \ , \ \ \ \lg u ,
1_\cA \rg \ =\ \varepsilon_\cU (u) 
\eena  
\\
The algebraic relations are then: 
\begin{eqnarray}
&& [\tilde D, C]  = -2C \ , \quad
   [B,C]  = \tilde D \ , \quad [B,C]_+  = \tD^2 \\
&& [\tilde D, B]  = 2 B {\tilde D}^2 \ ,
\quad [\tilde D, B]_+  = 0 \ , \quad   
   \tD^3  = \tD \ , \ \ \  C^2  = 0\ , \\ 
&& [\tA,Z]  = 0 \ , \quad Z=B,C,\tD \ , \\ 
&& E Z  = Z E  = 0 \ , \quad Z=\tA,B,C,\tD \qquad 
\end{eqnarray}
with co-algebraic structure:
\begin{eqnarray*}
  && \Delta(\tilde A)  = \tilde A\otimes 1_\cU + 1_\cU \otimes \tilde
  A, \quad  
  \Delta(B)  = B\otimes 1_\cU + 1_\cU\otimes B, \\
 && \Delta(C)  = C\otimes E + E\otimes C, \ \quad
 \Delta(\tilde D)  = \tilde D\otimes E + E\otimes \tilde D, \quad
 \Delta(E)  = E\otimes E 
 \\
 &&\eps_\cU (Z)  = 0 \ , \quad Z  = \tA, B, C, \tD \ ,\quad 
 \eps_\cU (E)  = 1 
\end{eqnarray*}

The extra operator $E$ is defined by: 
{~$\lg E,1_A\rg =1$}
 with all
other pairings being zero.
Strictly speaking the above algebra is in duality with the 
factor-algebra since it has zero pairings with the ideals. 
The full dual is infinitely generated and is under investigation.

Let us make a comparison with FRT duality \cite{FRT}. 
We use the relation 
$\lg\ L^\pm \ ,\ T \ \rg  = R^\pm $,
where:
$R^+\ \equiv\ P\,R\,P = R (-h)\ , \quad 
R^-\ \equiv\ R^{-1} $.
The comparison leads to
\begin{eqnarray}
 && L^{\pm}_{11}  = L^{\pm}_{22}  = e^{-h B} \ ,\quad
 L^{\pm}_{12}  = ((h_\pm + h^2)B + h \tilde A) e^{-h B} 
\end{eqnarray}
where $\ h_+  = h_3\ $ and $\ h_-  = -h_3 -2h^2$.
We see that $\tA$ and $B$ are taken into account in this formalism, 
but that it says nothing about  the generators $\ C,\tD\ $.

\subsection{$R$-matrix minimal polynomials and quantum planes}

In order to address the question of the quantum planes 
corresponding to the exotic bialgebras we have to know the 
minimal identity 
relations which the $R$-matrices fulfil. As we know the $R$-matrices 
producing deformations of the $GL(2)$ and $GL(1|1)$ fulfil second
order relations. However, in the cases at hand we have higher 
order relations. We have:
\eqna{hwrm} 
&& ( \hR - \id )^2 (\hR + \id) =  0 , \quad h_1 = -h_2 =
h,  h_3 \neq -h^2   \qquad\\
&& ( \hR - \id )^3 (\hR+ \id) =  0 , \quad h_1 + h_2 \neq
0 \\ 
&& ( \hR - \id ) (\hR + \id) =  0 , \quad h_1 = -h_2 =
h, h_3 = -h^2   \qquad
\eena 
where $\hR \equiv PR$, $R =R_{H2,3}\,$, 
\ $\id$\ is the $4\times 4$ unit matrix. 
The first case is E1, 
the second case is E2, 
the third case is 
the Jordanian subcase which produces  $GL_{h,h}(2)$.

\def\tB{\tilde B}
\def\tC{\tilde C}
\def\z{\zeta} 

To derive the corresponding quantum planes we shall apply the
formalism of \cite{WZ}.  The commutation relations 
between the coordinates $\ z^i\ $ and differentials $\ \z^i\ $, 
($i=1,2$), are given as follows:
\eqn{qpl} 
z^i z^j  = {\cP}_{ijk\ell}\, z^k z^\ell 
\ ,\quad
\zeta^i \zeta^j  = - {\cQ}_{ijk\ell}\, \zeta^k \zeta^\ell 
\ ,\quad
z^i \zeta^j  = {\cQ}_{ijk\ell}\, \zeta^k z^\ell 
\end{equation} 
where the operators $\ \cP,\, \cQ\ $ are functions of $\ \hR\ $
and must satisfy:
\eqn{orth} (\cP - \id)\,(\cQ + \id) = 0\ . \end{equation} 
Thus, there are different choices: four for E1, six for E2 
(and just one for $GL_{h,h}(2)$).
Choosing $\ \cP - \id = (\hR - \id)^a\
$ with $a=2,3,1$, respectively, and $\ \cQ = \hR\ $ in all cases.
and denoting $\ (x,y) = (z^1,z^2)\ $ we obtain 
\eqn{qpc} x y - y x = h y^2\ , \qquad h_1 = - h_2=h \ 
\end{equation} 
for E1 (and the Jordanian), or
{
\eqn{qpb} 
x y - y x = \frac12 (h_1-h_2) y^2\ , \qquad h_1 \neq - h_2 \ 
\end{equation} 
}
for E2. 
We note that the quantum planes corresponding to the three cases 
are not essentially different.

Denoting $\ (\xi,\eta) = (\zeta^1,\zeta^2)\ $ we obtain 
\eqn{dfb}
 \xi^2 + \frac{h_1 - h_2}{2}\ \xi \eta =0 \ ,\quad
\eta^2 =0 \ ,\quad
\xi \eta  = - \eta \xi 
\end{equation} 
for E2, while for E1 
$\xi^2 + h\, \xi \eta  = 0 $,
which is valid also for the Jordanian subcase. 

Finally, for the coordinates-differentials relations we obtain
\begin{alignat}{2}
& x \xi  = \xi x + h_1 \xi y + 
h_2 \eta x + h_3 \eta y \ ,\quad &&
x \eta  = \eta x + h_1 \eta y \ ,\\ 
& y \xi  = \xi y + h_2 \eta y \ ,\quad &&
y \eta  = \eta y 
\end{alignat}
\def\tb{\tilde b}\def\tc{\tilde c}
\def\half{{\textstyle{\frac{1}{2}}}}

\section{The S03 case}

\subsection{Algebraic relations}

The algebra relations in the S03 case are:
\begin{alignat}{3}
 \label{eq:S03trel}
 & \tb^2 = \tc^2 = 0 \ ,\qquad&&
 \ta\td = \td\ta = 0 \ ,\qquad&&
 \ta\tb = 0 \ , \\
 & \tb\td = 0\ ,\qquad&&
 \td\tc = 0\ , &&
 \tc\ta = 0\ . 
\end{alignat}
where: ~$ \tilde{b}  = \half (b+c)$, ~
$\tilde{c}  = \half (b-c)$. 

There is no PBW basis in this case. 
Indeed, the ordering  is cyclic: 
\eqn{ord} \ta ~>~ \tc ~>~ \td ~>~ \tb ~>~ \ta 
\end{equation} 
Thus, the basis consists of building blocks like ~$\ta^k\,
\tc\, \td^\ell\, \tb$~ and cyclic. 
Explicitly the basis can be
described by the following monomials:
\begin{subequations} \label{bas}
\begin{alignat*}{2}
&\ta^{k_1}\, \tc\, \td^{\ell_1}\, \tb\, \cdots\, 
 \ta^{k_n}\, \tc\, \td^{\ell_n}\, \tb\, \ta^{k_{n+1}} \ , \qquad
&&\td^{\ell_1}\, \tb\, \ta^{k_1}\, \tc\, \cdots\, 
 \td^{\ell_n}\, \tb\, \ta^{k_n}\ , \\ 
& \ta^{k_1}\, \tc\, \td^{\ell_1}\, \tb\, \cdots\, 
 \ta^{k_n}\, \tc\, \td^{\ell_n}\ ,
&&\td^{\ell_1}\, \tb\, \ta^{k_1}\, \tc\, \cdots\, 
 \td^{\ell_n}\, \tb\, \ta^{k_n}\, \tc\, \td^{\ell_{n+1}} \ ,
\end{alignat*}
\end{subequations}
where in all cases ~$n\,, k_i\,,\ell_i\, \in \, \ZZ_+\,$. 
\\
The algebra in duality is given by
\begin{eqnarray}
  \label{eq:s03}
  && [\tA,Z] = 0\ , \quad Z=\tB,\tC\ , \quad 
  \tA\tD = \tD\tA = \tD^3 =
  \tB^2 \tD = \tD \tB^2 = \tD \ , \nonumber\\
  && [\tB,\tC] = -2\tD \ ,\quad
  \tD\tB = - \tB\tD = \tC\tD^2 =\tD^2 \tC \ ,\quad
  \{\tC,\tD\} = 0\ , \nonumber\\
  & & \tB^2 + \tC^2 = 0\ ,\quad
  \tB^3 = \tB\ , \quad
  \tC^3 = - \tC\ , \quad
  \tB^2 \tA = \tA\ ,
\end{eqnarray}
with coproduct: 
\begin{eqnarray}
 \Delta_\cU(\tA) &=& \tA \otimes 1_\cU + 1_\cU\otimes \tA \\
 \Delta_\cU(\tB) &=& \tB \otimes 1_\cU + (1_\cU-\tB^2)\otimes \tB \\
 \Delta_\cU(\tC) &=& \tC \otimes (1_\cU-\tB^2) + 1_\cU\otimes \tC \\
 \Delta_\cU(\tD) &=& \tD \otimes (1_\cU-\tB^2) + (1_\cU-\tB^2)\otimes
 \tD \\ 
\eps_\cU (Z) &=& 0 \ , \qquad Z  = \tA, \tB, \tC, \tD \ .
\end{eqnarray}
$\tA$, ~$\tB^2=-\tC^2$~ and ~$\tD^2$~ are Casimir operators. 
\\
The bialgebra $s03$ is not a Hopf algebra (since there is no
antipode).

The algebra generated by the 
generator $\ \tA\ $ is a sub-bialgebra of ~$s03$. 
The algebra $s03'$ generated by $\tB,\tC,\tD$ 
is a nine-dimensional sub-bialgebra of ~$s03$~ with PBW basis:
\eqn{pbw:s03} 
1_\cU\,,\ \tB,\ \tC,\ \tD,\ \tB\tC,\ \tB\tD,\ \tD\tC,\ \tB^2,\
\tD^2 \end{equation} 
The algebra ~$s03$~ is not the direct sum of the
two subalgebras described above since both
subalgebras have nontrivial action on each other, e.g., 
~$\tB^2\tA ~=~ \tA$, ~$\tA\tD ~=~ \tD$. The algebra ~$s03$~ is a
nine-dimensional associative algebra over the central algebra
generated by $\tA$.
\\[2mm]
\noindent
Let us again make a comparison with FRT duality. 
\\
$L^\pm$ are matrices of operators $L^\pm_{ij}$ ($i,j=1,2$) 
satisfying the relations
\begin{eqnarray}
  \label{eq:RLL}
  R^+ L^+_1 L^+_2 = L^+_2 L^+_1 R^+ \\
  R^+ L^-_1 L^-_2 = L^-_2 L^-_1 R^+ \\
  R^+ L^+_1 L^-_2 = L^-_2 L^+_1 R^+ 
\end{eqnarray}
with $L_1\equiv L\otimes 1$, $L_2\equiv 1\otimes L$.
Explicitly, these RLL relations read
\begin{alignat}{2}
  \label{eq:RLLij}
  & (L^{\pm}_{11})^2 = (L^{\pm}_{22})^2 
  &\qquad\qquad
  & [L^{\pm}_{11}, L^{\pm}_{22} ] = 0 \nonumber\\
  & (L^{\pm}_{12})^2 = - (L^{\pm}_{21})^2 
  &
  & [L^{\pm}_{12}, L^{\pm}_{21} ]_+ = 0 \nonumber\\
  & L^{\pm}_{11} L^{\pm}_{12} = L^{\pm}_{22} L^{\pm}_{21} 
  &
  & L^{\pm}_{11} L^{\pm}_{21} = L^{\pm}_{22} L^{\pm}_{12} \nonumber\\
  & L^{\pm}_{12} L^{\pm}_{11} = - L^{\pm}_{21} L^{\pm}_{22} 
  &
  & L^{\pm}_{12} L^{\pm}_{22} = - L^{\pm}_{21} L^{\pm}_{11} 
\end{alignat}
and for the $RL^+ L^-$ ones
\begin{equation}
  \label{eq:RL+L-}
  L^+_{ij} L^-_{kl} - L^-_{ij} L^+_{kl}
  + \theta_i L^+_{\bar ij} L^-_{\bar kl}  + \theta_j L^-_{i\bar j}
  L^+_{k\bar l} = 0
\end{equation}
with $\bar n \equiv 3-n$, $\theta_1=1$, $\theta_2=-1$.
In the case of $s03$, the FRT relations in the dual algebra are richer than
the relations given using only (\ref{dub}). They will indeed lead to  more
irreducible finite dimensional representations.
\\
A convenient basis is given by
\begin{alignat}{2}
  &\tL^\pm_{11} = L^\pm_{11} + L^\pm_{22}
  &\qquad\qquad
  &\tL^\pm_{22} = L^\pm_{11} - L^\pm_{22}
  \nonumber\\
  &\tL^\pm_{12} = L^\pm_{12} + L^\pm_{21}
  &
  &\tL^\pm_{21} = L^\pm_{12} - L^\pm_{21}
\end{alignat}
In this basis, the relations (\ref{eq:RLLij}) read
\begin{alignat}{2}
  \label{eq:RtLtL1}
  &\tL^\pm_{11}  \,\tL^\pm_{22} = 0
  &\qquad\qquad\qquad
  &\tL^\pm_{22} \,\tL^\pm_{11} =0
  \nonumber\\
  &(\tL^\pm_{12})  ^2 = 0
  &
  &(\tL^\pm_{21}) ^2 = 0 
  \nonumber\\
  &\tL^\pm_{11}  \,\tL^\pm_{21} = 0
  &
  &\tL^\pm_{12} \,\tL^\pm_{11} =0
  \nonumber\\
  &\tL^\pm_{21}  \,\tL^\pm_{22} = 0
  &
  &\tL^\pm_{22} \,\tL^\pm_{12} =0
\end{alignat}
whereas the relations (\ref{eq:RL+L-}) become
\begin{alignat}{2}
  \label{eq:RtLtL2}
  &
  [\tL^+_{11}, \tL^-_{11}] = 0
  &\qquad\qquad
  &
  \tL^-_{21} \,\tL^+_{11} = \tL^+_{21} \,\tL^-_{11} 
  \nonumber\\
  &
  \tL^-_{11} \,\tL^+_{12} = \tL^+_{11} \,\tL^-_{12} 
  &\qquad\qquad
  &
  \tL^-_{21} \,\tL^+_{12} = \tL^+_{21} \,\tL^-_{12} 
  \nonumber\\
  &
  \tL^-_{11} \,\tL^+_{21} = \tL^+_{21} \,\tL^-_{22} 
  &\qquad\qquad
  &
  \tL^-_{21} \,\tL^+_{21} = - \tL^+_{11} \,\tL^-_{22} 
  \nonumber\\
  &
  \tL^-_{11} \,\tL^+_{22} = \tL^+_{21} \,\tL^-_{21} 
  &\qquad\qquad
  &
  \tL^-_{21} \,\tL^+_{22} = - \tL^+_{11} \,\tL^-_{21} 
  \nonumber\\[3mm]
  &
  \tL^-_{12} \,\tL^+_{11} = - \tL^+_{22} \,\tL^-_{12} 
  &\qquad\qquad
  &
  \tL^-_{22} \,\tL^+_{11} = \tL^+_{12} \,\tL^-_{12} 
  \nonumber\\
  &
  \tL^-_{12} \,\tL^+_{12} = - \tL^+_{22} \,\tL^-_{11} 
  &\qquad\qquad
  &
  \tL^-_{22} \,\tL^+_{12} = \tL^+_{12} \,\tL^-_{11} 
  \nonumber\\
  &
  \tL^-_{12} \,\tL^+_{21} = \tL^+_{12} \,\tL^-_{21} 
  &\qquad\qquad
  &
  \tL^-_{22} \,\tL^+_{21} = \tL^+_{22} \,\tL^-_{21} 
  \nonumber\\
  &
  \tL^-_{12} \,\tL^+_{22} = \tL^+_{12} \,\tL^-_{22} 
  &\qquad\qquad
  &
  [\tL^+_{22}, \tL^-_{22}] = 0
\end{alignat}

\noindent
Denote for $n \geq 1$:
\begin{eqnarray} 
  F_n(k_i;l_i) \equiv \prod^n_{i=1}\tilde{L}^{+
    k_i}_{11} \tilde{L}^+_{12} \tilde{L}^{+
    l_i}_{22}\tilde{L}^{+}_{21} \\ 
  G_n(l_i; k_i) \equiv \prod^n_{i=1}
  \tilde{L}^{+ l_i}_{22} \tilde{L}^+_{21} \tilde{L}^{+ k_i}_{11}
  \tilde{L}^+_{12} 
\end{eqnarray} 
and for $n=0$ 
\begin{equation} 
  F_0(k_i;l_i) \equiv 1; \ \ \ G_0(l_i;k_i) \equiv 1
\end{equation} 
The basis elements of the algebra generated by the $\tL^+$'s 
are (following \cite{ACDM2})
\begin{eqnarray} 
  \label{eq:basisB+}
  F_n(k_i;l_i)\tilde{L}^{+k_{n+1}}_{11}; \ \ \ F_{n-1}(k_i;l_i)
  \tilde{L}^{+k_n}_{11}\tilde{L}^+_{12}\tilde{L}^{+l_n}_{22}; 
  \nonumber\\
  G_n(l_i;k_i) \tilde{L}^{+l_{n+1}}_{22};  \ \ \
  G_{n-1}(l_i;k_i)\tilde{L}^{+l_n}_{22}\tilde{L}^+_{21}
  \tilde{L}^{+k_n}_{11}  
\end{eqnarray} 
Defining also $  K_n = \sum^n_{i=1} k_i, \ \ \ L_n = \sum^n_{i=1} l_i $
the actions of generators $\tL^-$ on the basis elements are, e.g.,
\begin{eqnarray} 
  \tilde{L}^-_{11} F_n(k_i;l_i)&=& F_{n-1}(k_1
  +1,k_i;l_i) \tilde{L}^{+ k_n}_{11} \tilde{L}^+_{12} \tilde{L}^{+
    l_n}_{22} \tilde{L}^-_{21}  \nonumber \\ \tilde{L}^-_{12}
  F_n(k_i;l_i) &=& (-1)^{K_n + L_n +1} G_{n-1}(k_1+1,k_i;l_i) 
  \tilde{L}^{+ k_n}_{22} \tilde{L}^+_{21}
  \tilde{L}^{+ l_n}_{11} \tilde{l}^-_{22} \nonumber\\ 
  \tilde{L}^-_{21} F_n(k_i;l_i) &=& G_n(0,..,l_{n-1};k_i) \tilde{L}^{+
    l_n}_{22} \tilde{L}^-_{21}  \nonumber \\ 
  \tilde{L}^-_{22}
  F_n(k_i,;l_i) &=& (-1)^{K_n + L_n}F_n(0,..,l_{n-1};k_i) \tilde{L}^{+
    l_n}_{11} \tilde{L}^-_{22}  
\end{eqnarray} 
\begin{eqnarray} 
  \tilde{L}^-_{11} G_n(l_i;k_i) &=&
  (-1)^{K_n+L_n}G_n(0, ..,k_{n-1};l_i) \tilde{L}^{+
    k_n}_{22}\tilde{L}^-_{11} \nonumber\\  
  \tilde{L}^-_{12}
    G_n(l_i;k_i) &=&  F_n(0,..,k_{n-1};l_i) \tilde{L}^{+
    k_n}_{11}\tilde{l}^-_{12}  
\end{eqnarray} 
\begin{eqnarray} 
  \tilde{L}^-_{21} G_n(l_i;k_i) &=& (-1)^{K_n+L_n+1}
  F_{n-1}(l_1+1,l_i;k_i) 
  \tilde{L}^{+ l_n}_{11} \tilde{L}^+_{12}
  \tilde{L}^{+k_n}_{22}\tilde{L}^-_{11} \nonumber \\ 
  \tilde{L}^-_{22} G_n(l_i;k_i) &=&
  G_{n-1}(l_1+1,l_i;k_i)\tilde{L}^{+
    l_n}_{22}\tilde{L}^+_{21}\tilde{L}^{+k_n}_{11} \tilde{L}^-_{12}
\end{eqnarray} 
These equations allow one to order the $\tL^-$ with respect to the
$\tL^+$. For the $\tL^-$ among themselves, there exists a basis
similar to (\ref{eq:basisB+}).

\subsection{Finite dimensional irreducible representations}

Let us first consider the algebra generated by $A$, $B$, $C$, $D$ and the 
relations (\ref{eq:s03}). 
Since $\tD^3=\tD$, there exists a weight vector ~$v_0$~ such that:
$ \tD\,v_0 ~=~ \lambda\,v_0 $
where ~$\lambda^3=\lambda$.
The finite  dimensional irreps are then
\begin{itemize}
\item one-dimensional trivial 
\item two-dimensional with Casimir values ~$\mu,1,1$~ 
for ~$\tA, \tB^2,\tD^2$, respectively, $\mu\in \CC$. 
\item one-dimensional with Casimir values ~$\mu,1,0$\\
for ~$\tA, \tB^2,\tD^2$, respectively, $\mu\in \CC$. 
\end{itemize}

A two-dimensional representation for the algebra generated by the
$L$-operators is provided by the $R$-matrix itself, setting $\pi(L^+)
= R_{21}$, $\pi(L^-) = R^{-1}$ 
(see \cite{FRT,AC0401})
\begin{alignat}{2}
  \label{eq:RRRlinear}
  &
  \pi(L^\pm_{11}) 
  =
  \left(
    \begin{array}{cc}
      1 &     0  \cr
      0 & \mp 1  
    \end{array}
  \right)   
  &\qquad\qquad
  &
  \pi(L^\pm_{12}) 
  =
  \left(
    \begin{array}{cc}
      0 & \pm 1  \cr
      1 &     0  
    \end{array}
  \right)   
  \nonumber\\
  &
  \pi(L^\pm_{21}) 
  =
  \left(
    \begin{array}{cc}
          0 & 1  \cr
      \mp 1 & 0  
    \end{array}
  \right)   
  &\qquad
  &
  \pi(L^\pm_{22}) 
  =
  \left(
    \begin{array}{cc}
      \pm 1 & 0  \cr
          0 & 1  
    \end{array}
  \right)   
\end{alignat}
Note that this does not exhaust the set of two-dimensional representations.

Let $N_1$ and $N_2$ be two  non negative integers.
Here is an example (in the $\tL$ basis) of a finite dimensional irreducible
representation of arbitrary dimension $N_1+N_2$. 
\begin{eqnarray}
  \label{eq:reprexample}
  &&
  \pi(\tL_{11}) =
  \diag(\rho_1,\cdots,\rho_{N_1},\underbrace{0,\cdots,0}_{N_2}) 
  \qquad\qquad \rho_i\neq\rho_j \qmbox{\rm for}  i\neq j\\ 
  &&\pi(\tL_{22}) =
  \diag(\underbrace{0,\cdots,0}_{N_1},\lambda_1,\cdots,\lambda_{N_2})
  \qquad\qquad \lambda_i\neq\lambda_j \qmbox{\rm for}  i\neq j
  \\ 
  &&\left( \pi(\tL_{12}) \right)_{ij} \neq 0 \qmbox{\rm iff} 
  i\in\{1,\cdots,N_1\}, \quad
  j\in\{N_1+1,\cdots,N_1+N_2\} \qquad
  \\
  &&\left( \pi(\tL_{21}) \right)_{ij} \neq 0 \qmbox{\rm iff} 
  i\in\{N_1+1,\cdots,N_1+N_2\}, \quad
  j\in\{1,\cdots,N_1\}
\end{eqnarray}

\subsection{Baxterisation}

We introduce the following Ansatz (choosing a convenient
normalisation): 
\begin{equation}
  \label{eq:am6}
  \hat{R}(x)=I+c(x)\hat{R}
\end{equation}
and we try to find $c(x)$ such that $\hat{R}(x)$ would satisfy 
the parametrised Yang--Baxter equation:
\begin{equation}
  \hat{R}_{(12)}(x)\hat{R}_{(23)}(xy)\hat{R}_{(12)}(y) =
  \hat{R}_{(23)}(y)\hat{R}_{(12)}(xy)\hat{R}_{(23)}(x)  
\end{equation}

The result is:
\begin{equation}
  \label{eq:am11}
  \hat{R}(x)= (\sqrt{2}x)^{-1} \hat{R} +(\sqrt{2}x)\hat{R}^{-1} 
  = 
  \frac{1}{\sqrt{2x}}
  \left(
    \begin{array}{cccc}
      x+1 & 0 & 0 & 1-x \cr
      0 & x+1 & x-1 & 0 \cr
      0 & 1-x & x+1 & 0 \cr
      x-1 & 0 & 0 & x+1                                                   
    \end{array}
  \right)     
\end{equation}

\subsection{Spectral decomposition and noncommutative planes}

The  minimal polynomial identity satisfied by $\cR$ is
\begin{equation}
  \label{eq:am1}
 \hat{R}^2 - 2\hat{R} +2I = 0 
\end{equation} 
and we have the spectral decomposition: 
\begin{equation}
  \label{eq:am5}
  \hat{R}=(1-i)P_{(+)}+(1+i)P_{(-)}  = (1+i)I-2iP_{(+)} 
\end{equation}
\begin{equation}
  \label{eq:am3}
  \mbox{where}\qquad
  P_{(\pm)} ~\equiv ~ \frac{1}{2} (I \pm i (\hat{R} -I)) 
\end{equation}
are projectors resolving the identity:
$P_{(i)}P_{(j)} = {{\delta}_{ij}}P_{(i)},\quad i,j = \pm $,
$  P_{(+)}+P_{(-)}=I $.
The quantum plane relations in the case of S03 are then
\begin{eqnarray}
  \label{eq:am28}
  {x_{1}}^2=x_{1}x_{2}, \quad {x_{2}}^2= -x_{2}x_{1}
  \ ,\quad
  {{\xi}_{1}}^2=-{\xi}_{1}{\xi}_{2}, \quad {{\xi}_{2}}^2= {\xi}_{2}{\xi}_{1}
  \\
  x_{1}{\xi}_{1}=(\nu-1){\xi}_{1}x_{1} + \nu{\xi}_{1}x_{2}
  \ ,\quad
  \label{eq:am30}
  x_{1}{\xi}_{2}=(\nu-1){\xi}_{1}x_{2} + \nu{\xi}_{1}x_{1}
  \\
  x_{2}{\xi}_{1}=(\nu-1){\xi}_{2}x_{1} - \nu{\xi}_{2}x_{2}
  \ ,\quad
  x_{2}{\xi}_{2}=(\nu-1){\xi}_{2}x_{2} - \nu{\xi}_{2}x_{1}
\end{eqnarray}
where $\nu$ is arbitrary real parameter.

\section{The S14  case}

The relations in the S14 case are:
\begin{align*}
 \label{eq:S14trel}
 & \tb\tc + \tc\tb = 0 \qquad\qquad
 \ta\td + \td\ta = 0 \nonumber\\
 & \ta\tb = \tb\ta = \ta\tc = \tc\ta = \tb\td = \td\tb = 
 \tc\td = \td\tc = 0  
\end{align*}
and the dual algebra is:
\eqna{eq:dualS14}
&& \tC = \tD\tB = -\tB\tD \ ,\quad 
[\tA,\tD] = 0 \ ,\quad
E Z  = Z E  = 0
\\
&& \tA\tB = \tB\tA = \tD^2 \tB = \tB^3 = \tB  \ , \quad Z=\tA,\tB,\tD \ . 
\eena 
The dual coalgebra is given by (with $K \equiv (-1)^{\tA}$)
\eqna{eq:dualco14}
 \Delta_\cU(\tA) &=& \tA \otimes 1_\cU + 1_\cU\otimes \tA \ ,\quad
 \Delta_\cU(\tB) = \tB \otimes E + E\otimes \tB \\
 \Delta_\cU(\tD) &=& \tD \otimes K + 1_\cU\otimes \tD \ , \qquad 
 \Delta(E) = E\otimes E 
 \\
 \eps_\cU (Z) &=& 0 \ \ \mbox{for} \quad Z  = \tA, \tB, \tD \ ; \qquad
 \eps_\cU (E) = 1
 \eena 
The irreducible representations of ~$s14$~ 
follow the classification
\begin{itemize}
\item one-dimensional with Casimir values ~$\mu,0,\lambda^2$~ 
for ~$\tA, \tB^2,\tD^2$, respectively, $\mu,\lambda\in \CC$. 
\item two two-dimensional with all Casimirs 
~$\tA, \tB^2,\tD^2$~ having the value~1. 
\end{itemize}

\section*{Acknowledgments}
This work was supported in part by
the CNRS-BAS\ France/Bulgaria agreement number 12561. 
It was also inspired by the TMR Network
EUCLID: ``Integrable models and applications: from strings to
condensed matter'', contract number HPRN-CT-2002-00325.



\begin{thebibliography}{10}

\bibitem{Hlavaty}
  L. Hlavat\'y, 
  {\it J. Phys.} {\bf A20}, 1661 (1987); 
  {\it J. Phys.} {\bf A25}, L63 (1992).

\bibitem{Hietarinta}
 J. Hietarinta, 
 {\it J. Math. Phys.} {\bf 34}, 1725 (1993). 

\bibitem{ACDM1} D. Arnaudon, A. Chakrabarti, V.K. Dobrev and S.G.
  Mihov, 
  {\it J. Phys.} {\bf A34},
  4065  (2001); \texttt{math.QA/0101160}.
  
\bibitem{ACDM2} D. Arnaudon, A. Chakrabarti, V.K. Dobrev and S.G.
  Mihov, 
  {\it J. Math. Phys.} {\bf 43}, 6238  (2002); \texttt{math.QA/0206053}.
  
\bibitem{ACDM3}
  D. Arnaudon, A. Chakrabarti, V.K. Dobrev and S.G. Mihov,
  {\it Int. J. Mod. Phys.} {\bf A18}, 4201  (2003);
  \texttt{math.QA/0209321}.

\bibitem{WZ} J. Wess and B. Zumino, \textsl{Covariant differential
    calculus on the quantum hyperplane}, {\it Nucl. Phys. (Proc. Suppl.)}
  {\bf 18}, 302 (1990). 

\bibitem{FRT} L.D. Faddeev, N.Yu. Reshetikhin and L.A. Takhtajan,
  ``Quantization of Lie groups and Lie algebras'', 
{\it Alg. Anal.} {\bf 1}, 178-206 (1989)  (in Russian) and in: {\it
    Algebraic Analysis}, Vol. 1 (Academic Press, 1988) pp.
  129-139.

\bibitem{AC0401} A. Chakrabarti, ``A nested sequence of
    projectors and corresponding braid matrices $\hat R(\theta)$ (1)
    odd dimensions'',   \texttt{math.QA/0401207}. 

\end{thebibliography}
\end{document}